\documentclass{amsart}

\usepackage{amsmath,amssymb,url,algorithm2e}
\usepackage{graphicx}

\newcommand\Z{{\mathbb Z}}

\newcommand{\rnk}[2]{\lower 0.2ex\hbox{$#1$}\kern -.1em{\setminus}
                     \kern -.1em\raise 0.2ex\hbox{$#2$}}
\newcommand{\dnk}[3]{\lower 0.2ex\hbox{$#1$}\kern -.1em{\setminus}
                     \kern -.1em\raise 0.2ex\hbox{$#2$}
                     \kern -.1em{/}\kern -.1em\lower 0.2ex\hbox{$#3$}}
\newcommand{\lnk}[2]{\raise 0.2ex\hbox{$#1$}
                     \kern -.1em{/}\kern -.1em\lower 0.2ex\hbox{$#2$}}

\makeatletter
\def\mkdirsemi#1#2#3#4{\mathbin{\times\hbox to #4{}\hbox{\vrule
  \@width #1\@height #2\@depth #3}\hbox to .2em{} }}
\makeatother

\newcommand\SL{\mbox{\rm SL}}
\newcommand\Sp{\mbox{\rm Sp}}

\newcommand\gen[1]{\langle#1\rangle}
\newcommand\sz[1]{\left|#1\right|}
\newcommand\gesy[1]{\underline{\mathbf{#1}}}
\newcommand\mvec[1]{\underline{\mathbf{#1}}}

\begin{document}
\title{Constructive Membership Tests in Some Infinite Matrix Groups}
\author{Alexander Hulpke}
\address{Colorado State University,
Department of Mathematics,
1874 Campus Delivery,
Fort~Collins,
Colorado,
80523-1874
USA}
\email{hulpke@colostate.edu}

\begin{abstract}
We describe algorithms and heuristics that allow us to express arbitrary
elements of $\SL_n(\Z)$ and $\Sp_{2n}(\Z)$ as products of generators in
particular ``standard'' generating sets. For elements obtained
experimentally as random products, it produces product expressions whose
lengths are competitive with the input lengths.
This is the author's copy of~\url{https://doi.org/10.1145/3208976.3208983}
\end{abstract}

\maketitle

\section{Introduction}

The constructive membership problem, that is expressing an element $g$ of a
group $G$ as a word in a generating set (which might be user-chosen)
is one of the fundamental tasks of computational group theory. 
We call such a word a {\em factorization} of $g$ (with respect to the chosen
generating set).
In the case
of elementary abelian groups it is simply the well-studied problem of
solving a system of linear equations.

Another special case is discrete logarithm, that is expressing an element in
a cyclic group as a power of a chosen generator. This is a known, difficult
problem~\cite{schoofdiscrete}, thus
the best we can hope for are good heuristics rather than a general solution.

We note here that~\cite{babaibealsseress09} shows that for a large class of
groups, including matrix groups, discrete logarithm to be the only obstacle to
efficient group order and membership test calculations.
\smallskip

The application to puzzles~\cite{egnerpueschel} arguably has the
largest visibility for the general public.

In general, factorization underlies much of the functionality for group
homomorphisms~\cite{lgpraesoi} and is thus at the heart of many group theoretic
calculations.

While for problems such as the Cayley graph for Rubik's
cube~\cite{rokickicube,kunklecooperman} a shortest word expression is
the inherent aim, in most applications the goal is rather to
obtain a word expression that is ``reasonably short'' for practical purposes,
but without any guaranteed bound in relation to the optimal length
(or even just a straight line program).
\smallskip

For permutation groups, stabilizer chains~\cite{sims70}
provide a tool for obtaining such expressions~\cite{Minkwitz}. For finite
matrix groups, composition trees and constructive
recognition~\cite{baarnhielmholtleedhamobrien} are tools.
\medskip

The groups we are interested in here will be particular infinite matrix
groups over rings of integers,  namely $G=\SL_n(\Z)$ or $G=\Sp_{2n}(\Z)$
with particular generating sets.  This is motivated by recent
work~\cite{detinkoflanneryhulpke17} on finitely generated subgroups of these
groups:
Given a subgroup $S\le G$ given by generating
matrices, one often  would like to
determine whether $S$ has finite index in $G$, in which case $S$ is called {\em
arithmetic}. Calculation in finite images of $G$ allow
us~\cite{detinkoflanneryhulpke17} to
determine the index if it is known to be finite.
\smallskip

Determining {\em whether} the index is finite, however requires us\footnote{%
Structural arguments based on the existence of free subgroups show that
there cannot be deterministic, bounded-time finite index test. Any method
that has a chance of determining the index thus needs to share
characteristics of methods for subgroups of finitely presented groups.} to
verify the index in a finitely presented version of $G$ and thus poses the
task to express the generators of $S$ as words in a particularly chosen
generating set of $G$. While methods for word expression exist already for
the case of $\SL$, the author has been unable to find such methods for
$\Sp$ in the literature. The Section~\ref{xample} below will give an
example
(taken from~\cite{hofmanvanstraten})
of doing this using the approach presented in this paper.
\smallskip

While these groups clearly exist in arbitrary dimension, the questions and
concrete examples studied so far have been of rather limited dimension ($\le
8$). One reason for this is that products of elements of infinite matrix
groups usually very quickly produce large coefficients, and matrix
arithmetic itself becomes a bottleneck.

This paper thus is focusing on practically useful methods for small
dimensions, even if they scale badly for larger $n$.
\smallskip

This use of the factorization also indicates that the appropriate measure of
success is the length of the resulting words, rather than the time required
to obtain such a factorization: The time for the overall
calculation will be dominated by the coset enumeration, and shorter words
often make success of such an enumeration more likely.
\smallskip

We shall present algorithms that in experiments perform well under this
measure, though we cannot give a provable statement about the quality of the
word expression obtained.
In the case of $\Sp$, furthermore we shall present a
heuristic that has worked well for all examples tried, though we cannot
prove this statement in general.

\section{Two basic algorithms}

We start by fixing notation: We have a group $G$
with a generating
sequence $\mvec{g}=(g_1,\ldots,g_k)$. The task is to express an arbitrary $e\in
G$ as a {\em word} in $\mvec{g}$, that is a product of the elements in
$\mvec{g}$ and their inverses that equals $e$. We shall call such a word a
{\em word expression} for $e$.
The smallest number of factors possible in such a word expression for $e$ is
called the {\em word length} of $e$ (with respect to the generating set
$\mvec{g}$),
and such a word is called a {\em shortest word} for $e$.

To simplify notation, we shall also assume now that $\mvec{g}=\mvec{g}^{-1}$ is
closed under taking inverses.
\medskip

The {\em Cayley graph} $\Gamma$ of $G=\gen{\mvec{g}}$ is a digraph with vertex set
$G$ and, for $x,y\in G$, an edge $(x,y)$ labeled by $g$, existent iff
$xg=y$. The question for a word of minimal length expressing $e\in G$ thus
is the same as that of finding an (undirected) path in $\Gamma$ of shortest
length from $1_G$ to $g$.

Standard ``shortest path'' algorithms for graphs, such as~\cite{dijkstra}
then motivate an exhaustive search that ``floods'' the Cayley graph vertex
by vertex, starting with the identity and stopping once group element 
$e$  has been reached. The corresponding algorithm for word expression has been known
for a long time and is given as algorithm~\ref{brute}:
\smallskip

We use the notation $L[a]$ to get a list element associated to a group
element $a$, this will be implemented though appropriate data structures,
such as hashing.
\begin{algorithm}
\caption{Floodsearch}
\label{brute}
\SetKwInOut{Input}{Input}
\SetKwInOut{Output}{Output}
\Input{A group $G$ with generating set $\mvec{g}$ and $e\in G$}
\Output{A word expression in $\mvec{g}$ for $e$ or a memory overflow error}
Initialize $A:=\{(1_G)\}$\;
$P:=[]$, $P[1_G]$:=false
\tcp*{Marker whether an element was processed}
$W:=[]$; $W[1_G]:=\emptyset$
\tcp*{Word expressions for elements}
\If{$e=1$}{\Return{$\emptyset$}\;}
\While{Memory is not exhausted}{
Let $A'=\{a\in A\mid P[a]=\mbox{false}\}$\;
  \ForEach{$a\in A'$}{
    \ForEach{$x\in\mvec{g}$}{
      \If{$ax=e$}{
        \Return{$(W[a],x)$}\tcp*{Concatenate words}
      }
      \ElseIf{$ax\not\in A$}{
        add $ax$ to $A$\; 
	set $W[ax]:=(W[a],x)$\tcp*{Concatenate words}
	$P[ax]$:=false\;
      }
    }
    Set $P[a]$:=true\;
  }
}

\tcp{Stage 2: Word products}

\ForEach{$a\in A$}{
  \If{$ea\in A$}{
    \Return{$(W[ea],W[a^{-1}])$}\;
  }
}

\Return{Memory exhaustion failure}\;
\end{algorithm}

The first stage of this approach can also be considered as an orbit algorithm~\cite{holtbook},
calculating the (partial) orbit $A$ of $1_G$ under right multiplication by $G$.
In this form it is easily seen that that it is sufficient not to store full
word expressions $W$, but only the generator labeling the last edge of the
shortest path.
(In fact, following~\cite{coopermanfinkelstein92}, one can reduce the storage
requirement to 2 bits per element by indicating the length of the path
modulo 3.)
\smallskip

Fundamentally, this is one of the the only two known approaches that can
guarantee\footnote{E.g.  the calculation of the diameter of Rubik's
cube~\cite{rokickicube} ultimately builds on this algorithm} to find a
shortest word.  The other method would be to use a finite, length-based
confluent rewriting system for $G$. (In general we do not have good
confluent rewriting systems for arbitrary finite groups, furthermore in the
infinite case it is not even known whether such finite systems exist.)
\medskip

If memory is exhausted, we then can (this is {\tt Stage 2})
use the fact that all elements are invertible and that the Cayley
graph looks the same from every vertex to extend the radius by a factor two,
before failure:
Test whether the ball $A$ (around $1_G$) and the ball $eA$ (around $e$)
intersect:

If a word of shortest length is desired, we may not stop at the first word
that is found, but must run systematically through all combinations (or run
through pairs according to the length of the product).
\smallskip

The storage requirements, which are ${\mathcal O}(\sz{A})$, show that for
every group there is a maximal word length that can be tested for. Thus this
method can cater only for a finite number of elements in an infinite group.
Its use is rather are  as ``quality control'' of the produced word length
for other algorithms, or to find explicit word lengths for particular
elements.

\subsection{Modular reduction}
Another algorithm is specific to integral matrix groups:
Given $e\in\SL_n(\Z)$,
we find a word expression in a finite congruence image
$\SL_n(\Z/p\Z)$, for example using stabilizer chain methods.
Having found such an expression, we then
check whether this expression also holds in
characteristic zero. Otherwise we consider larger congruence images.
The implicit expectation here is that for sufficiently
large modulus $p$ no modular reduction happens in evaluating a word
expression for $e$ and the calculation modulo $p$ is in fact the same as the
calculation in $\Z$.
This is decried as algorithm~\ref{congrimg}.

\begin{algorithm}
\caption{Word by congruence image}
\label{congrimg}
\SetKwInOut{Input}{Input}
\SetKwInOut{Output}{Output}
\Input{A group $G\le\SL_n(\Z)$ with a generating set $\mvec{g}$ and $e\in G$}
\Output{A word expression in $\mvec{g}$ for $e$ or failure}
Let $p=3$\;
\While{$\sz{\SL_n(p)}$ is not too large}{
  Let $\varphi\colon G\to\SL_n(p)$ the congruence homomorphism\;
  Let $H=\gen{\varphi(\mvec{g})}$\;
  Let $w$ be a word expression for $\varphi(e)$ as a word in
  $\varphi(\mvec{g})$\;
  \If{$w$ evaluated in $\mvec{g}$ equals $e$}{
    \Return{$w$}\;
  }
  Increment $p$ to the next prime\;
}
\Return{failure}\;
\end{algorithm}

Despite its simplicity, this method often works well to find short word
expressions (as done for example in~\cite{detinkoflanneryhulpke14}
to find candidates
for generic word
expression that then are explicitly proven) and sometimes works faster (and for longer
word lengths) than the previous one. However there is no practical way to
determine a priori a small modulus that would guarantee success.

\section{Norm based methods for $\SL$}
\label{dosl}

We now consider the special case of $G=\SL_n(\Z)$ with generators being
elementary matrices.  We denote by $t_{i,j}$ the matrix that is the identity
with an extra entry one in position $i,j$ and set 
\[
\gesy{g}=\left\{t_{i,j}^{\pm 1}\mid q\le i\not=j\le n \right\}.
\]
Then~\cite{hahnomeara} (which is a basic linear algebra argument)
shows that $\SL_n(\Z)=\gen{\mvec{g}}$.

A word in these generators can be considered as performing a sequence of
elementary matrix operations, and the inverse of a word expression for $g\in
G$ would be a sequence of elementary operations that transform $g$ to the
identity, which is also its Hermite Normal Form (HNF).
(As the normal form is the identity, the calculation of Hermite Normal Form
is effectively the same as that of the Smith Normal Form in this case.)
\smallskip

Calculating Normal Forms of matrices is a classical problem in Computer
Algebra~\cite{storjohann98,saunderswan04}. If we perform such a calculation
and accumulate the sequence of
elementary operations not in transforming matrices, but as words, we
obtain a word expression in terms of elementary matrices.
We shall call this algorithm~\ref{hnfbased} the HNF-based algorithm.

\begin{algorithm}
\caption{The HNF-based algorithm}
\label{hnfbased}
\SetKwInOut{Input}{Input}
\SetKwInOut{Output}{Output}
\Input{A group $G\le\SL_n(\Z)$ with a generating set $\mvec{g}$ consisting
of elementary matrices, and $e\in G$}
\Output{A word expression in $\mvec{g}$ for $e$}
Calculate the HNF for $e$ and the transforming matrix $T$ such that
$Te=1$. While doing so keep $T$ as a word expression in the elementary
matrices $\mvec{g}$.\;
\Return{$T^{-1}$}
\tcp*{Use $T^{-1}$ since $T$ converts $e$ to $1$}
\end{algorithm}

An implementation of this algorithm was built on top of the {\sf GAP} \cite{GAP4}
implementation of
Hermite Normal Form. For matrices with moderate entries (respectively those
who have word length in the generators of not more than 20-30) it produces
satisfactory results, but not if examples of longer word length are
considered. This is because longer products correlate with larger
coefficients. For such matrices the first steps in a normal form calculation
are to reduce a row by subtracting the $k$-th multiple of another row,
typically for a large $k$. Such steps produce an elementary matrix in $k$-th power
and thus makes for very long words. This is corroborated by the
examples in section~\ref{xample}.
\smallskip

Note also that this approach only applies if the group is
generated by all elementary matrices. It thus is only applicable for $\SL$,
not subgroups thereof.
\medskip

We thus consider further strategies used for
calculating normal forms, rather than to utilize the forms themselves.

The starting observation is that matrix multiplication in characteristic
zero tends to produce a product that has larger entries than either factor.
Reducing the overall size of entries of the matrix thus is expected to be
more promising than trying to zero out off-diagonal entries systematically
row-by-row and column-by-column.
\smallskip

We shall use the (squared) 2-matrix norm $\|M\|^2=\sum_{i,j}
m_{i,j}^2$. A smaller norm corresponds to overall
smaller entries.
In fact, as we know the normal form to be the identity, we 
use the measure $\|M-I\|^2$ ($I$ being the identity matrix) in place of
$\|M\|^2$.

We shall denote $\|M-I\|^2$ from now on as {\em height}.

The algorithm for factorization now iterates a reduction process for the
entries of a matrix $a\in\SL_n(\Z)$, as given by algorithm~\ref{heightbase}:
We try to reduce matrix height by forming products with generators. In a
greedy algorithm we form products with all generators and choose the one
that produces the largest height reduction. If no such generator exists we
fall back on the proven HNF-based method.
\begin{algorithm}
\caption{Height-based reduction}
\label{heightbase}
\SetKwInOut{Input}{Input}
\SetKwInOut{Output}{Output}
\Input{A group $G\le\SL_n(\Z)$ with generating set $\mvec{g}$ (that is
assumed, but not required, to contain elementary matrices) and $e\in G$}
\Output{A word expression in $\mvec{g}$ for $e$ or failure}
Let $w=\emptyset$, $a:=e$\;

\While{$a\not=I$}{
  \ForEach{$g_i\in\mvec{g}$}{
    Calculate $\|a\cdot g_i-I\|^2$ and $\|g_i\cdot a-I\|^2$, and find for
    which $g_i$ and product order
    the value $m$ is minimal\;
  }
  \If{$m\ge\|a\|^2$}{
    Factor $a$ with the HNF-based method (algorithm~\ref{hnfbased}),
    obtaining a word $v$ for $a$\;
    \Return{$v\cdot w^{-1}$}\;
  }
  Replace $a$ with the product that produced minimal height\;
  Replace $w$ by (the corresponding) $(w,g_i)$, respectively $(g_i,w)$\;
}
\Return $w-1$\;
\end{algorithm}

Applying this algorithm to random elements of $\SL_n(\Z)$ produces in most
cases a significant reduction in the matrix coefficients, but do not reach
the identity before having to
default to algorithm~\ref{hnfbased}:

For example, let $\mvec{g}$ the set of all $4\times 4$ elementary matrices
and consider the element 
\[
a=\left(\begin{array}{rrrr}%
1&0&1&-1\\%
1&0&0&0\\%
0&-1&2&0\\%
0&-1&0&1\\%
\end{array}\right).
\]
Then $\|a-I\|=7$, but for no generator $g_i\in\mvec{g}$ we get a
smaller height.
\smallskip

In such a situation the matrix $a$ typically will have undergone prior
reduction and thus have comparatively small coefficients. Thus using the
HNF-based algorithm as fall-back is less likely to incur the word length penalty we noted
before.
\smallskip

Investigating this example further, we find (using algorithm~\ref{brute})
that $a$ can be written as a
product
$t_{1,4}^{-1}t_{2,1}t_{3,2}^{-1}t_{1,2}^{-1}t_{2,3}^{-1}t_{4,2}^{-1}$
of length $6$. If we calculate the heights of partial products of this word
we get $1,2,3,5,7,7$ if we take subwords starting from the left,
respectively $1,2,4,6,7,7$ from the right. The reason for
the failure of the height-based approach thus is that the first, as well as the
last factor of the product, do not increase the height.
\smallskip

This failure to reduce is the result of small height values and looking at
only single generators. If instead we would have considered products of
length $2$, we would have noticed a jump from $7$ to $\le6$ by multiplying
with a product of length $2$.
\smallskip

Obviously, one could try also longer products. In an ad-hoc compromise
between length and number of products we decided to consider
products of length up to $3$, as this includes conjugates of generators by
other generators.

Let 
\[
\mvec{h}=\mvec{g}\cup\mvec{g}^2\cup\mvec{g}^3.
\]
When the height-based reduction then reaches the stage at which no element of
$\mvec{g}$ reduces, we repeat
the same attempt of height reduction through generators, albeit with
$\mvec{h}$ in place of $\mvec{g}$.
To avoid a careless accumulation of longer products, we furthermore
weigh the height change achieved by the length of the product expression
used.

If use of the generating set $\mvec{h}$ achieves a height reduction we
change $a$ and $w$ accordingly. 
If also use of $\mvec{h}$ achieved no improvement we pass to the HNF-based
algorithm~\ref{hnfbased}. Otherwise
the calculation then continues again with
reductions by $\mvec{g}$. 
\medskip

In experiments with random input (see section~\ref{xample}) we found that the
words produced by this approach seemed to be of acceptable length.

The time taken in the examples considered was short enough that we did not
look into ways to speed up the calculation, though there are many obvious
ways to do so, e.g. by looking at changes locally rather than always
processing a whole matrix.

\section{The Symplectic Group}
\label{spsec}

The symplectic group of degree $2n$ is the group of matrices in
$\SL_{2n}(\Z)$ that preserve the bilinear form
\[
J=\left(\begin{array}{cc}0&I_n\\
-I_n&0\\
\end{array}
\right)
\]
with $I_n$ denoting an $n\times n$ identity matrix. Thus
\[
\Sp_{2n}(\Z)=\left(M\in\SL_{2n}(\Z)\mid MJM^T=J
\right).
\]

A presentation for $\Sp_{2n}(\Z)$  has been calculated by
Birman~\cite{birman71}, based on prior work by Klingen~\cite{klingen61} and
unpublished thesis work of Gold~\cite{goldthesis61}.
Klingen's results shows that the given elements
indeed generate $\Sp$, but is very much non-constructive. It thus
does not facilitate an algorithm for
decomposition in these generators.

The work in~\cite{birman71} minimally adjusts the generating set
of~\cite{klingen61} and uses
\[
\Sp_{2n}(\Z)=\gen{ Y_i, U_i, Z_j\mid 1\le i\le n, 1\le j\le n-1 }
\]
with $Y_i=t_{i,n+i}^{-1}$, $U_i=t_{n+i,i}$ and 
\[
Z_i=\left(t_{i+1,n+i}/t_{i+1,n+i+1}\right)^{t_{i,i+1}}
=\left(\begin{array}{cc}
I_n&B_i\\
0&I_n\\
\end{array}
\right)
\]
with $B_i$ the matrix with submatrix
$\left(\begin{array}{rr}-1&1\\
1&-1\\
\end{array}
\right)$ at positions $i,i+1$ along the diagonal (and all other entries
zero).
\medskip

The generators $Z_i$ are not elementary, which does not augur well for
simply
replicating the approach used for $\SL$. We thus note that by~\cite{hahnomeara} we
can generate $\Sp$ as well from a generating set consisting of short
products of elementary matrices, resulting in a second generating set that
overlaps with the previous one:
\begin{eqnarray*}
\Sp_{2n}(\Z)&=&\langle\left\{t_{i,n+j}t_{j,n+i}, t_{n+i,j}t_{n+j,i}
\mid 1\le i<j\le s \right\}\\
&&\cup\left\{ t_{i,n+i},t_{n+i,i}\mid 1\le i\le n \right\}\rangle.
\end{eqnarray*}
We however do not have a presentation in this second generating set (though
one could produce one through a modified Todd-Coxeter algorithm, albeit at
the cost of relator lengths).
We simply add these elements as further generators (together with relations
that express them as products in the original generating set).
\smallskip

To find the necessary product expressions,
we need to express the elements $t_{i,n+j}t_{j,n+i}$ (whose
factors lie outside $\Sp$) as product of our chosen (primary)
generators for $\Sp$.

In \cite{birman71} we already find an expression for some of these products:
For $i\le n-1$ we have that
\begin{eqnarray*}
t_{i,i+1}t_{n+i+1,n+i}^{-1}&=&Y_i^{-1}Y_{i+1}^{-1}U_{i+1}^{-1}Y_{i+1}^{-1}Z_{i}U_{i+1}Y_{i+1}\qquad\mbox{and}\\
t_{i+1,i}t_{n+i,n+i+1}^{-1}&=&Y_{i+1}Y_iU_iY_iZ_i^{-1}U_i^{-1}Y_i^{-1}
\end{eqnarray*}
(Algorithm~\ref{brute} confirms that these are word
expressions of minimal length.)

We similarly used algorithm~\ref{brute} to suggest short expressions for other
products, and obtained for $3\le i\le n$ that:
\begin{eqnarray*}
t_{i-2,i}t_{n+i-2,n+i}^{-1}&=&
[Y_{i-1}Z_{i-1}^{-1}U_iY_i,U_{i-1}Y_{i-1}Z_{i-2}^{-1}U_{i-1}^{-1}],\\
\mbox{and}&&\\
t_{i,i-2}t_{n+i,n+i-2}^{-1}&=&
[Y_{i-1}Z_{i-2}^{-1}U_{i-2}Y_{i-2},\\
&&\qquad U_{i-1}Y_{i-1}Z_{i-1}^{-1}U_{i-1}^{-1}]
\end{eqnarray*}
with $[a,b]=a^{-1}b^{-1}ab$ denoting the commutator.

These expressions are easily verified in general by considering the images
of standard basis vectors under the left hand size products and the right
hand side products.

The identity $t_{j,i+j}=[t_{j,i+j-1},t_{i+j-1,i+j}]$ finally allows us to
form all other products $t_{i,j}t_{n+j,n+i}^{-1}$ as commutators of
products with smaller index difference.
\smallskip

From now on $\mvec{\tilde g}$ shall denote this extended generating set,
consisting of the $U_i$, $Y_i$, $Z_j$ and products
$t_{i,n+j}t_{j,n+i}, t_{n+i,j}t_{n+j,i}$ (and inverses thereof).

We experimented with algorithm~\ref{heightbase} with this generating set
$\mvec{\tilde g}$ (of course with the HNF-based
method replaced by an error message) on a number of random elements.

The results of these experiments were disappointing: Almost all elements we
tried reduced only partially and still left matrices with large entries for
which no further reduction process could be found, not even by introducing
further short products. 

\subsection{Decomposing the Symplectic group}

Thus a more more guided reduction, adapted to the structure of
the symplectic group, is required. We should emphasize however that the
following approach is purely heuristic in that we have is no proof of it
succeeding in general. In a large number of examples in small dimension, we
however also were unable to find a single example in which the approach
failed.
\medskip

We start with a structural observation:
The products $t_{i,j}t_{n+j,n+i}^{-1}$ for all $1\le i\not=j\le n$
clearly generate a subgroup
\[
S=\left\{\left(\begin{array}{cc}M&0\\
0&M^{-1}\end{array}\right)\mid M\in\SL_n(\Z)
\right\}\le\Sp_{2n}(\Z).
\]
such that $S\cong\SL_n(\Z)$. The definition of the symplectic group
(and the fact that $\Z$ has only two units) shows that for
\[
T=\left\{\left(\begin{array}{cc}\star&0\\
0&\star\\
\end{array}\right)\in\SL_{2n}(\Z)
\right\}
\]
and 
\[
R=T\cap\Sp_{2n}(\Z)=
\left\{\left(\begin{array}{cc}M&0\\
0&M^{-1}\\
\end{array}\right)\in\SL_{2n}(\Z)
\right\}
\]
we have that $S\le R$ is of index $2$.
\smallskip

This subgroup $R$ lies at the heart of the new approach.
If we have an element $e\in R$, we can use multiplication by
\[
(Y_1^2U_1)^2=\left(\begin{array}{cc}A&0\\
0&A\end{array}\right)\mbox{\ with\ }
A=\left(\begin{array}{cccc}
-1&&0&\\
&1&&\\
0&&\ddots&\\
&&&1\\
\end{array}
\right)
\]
to obtain $e'\in S$. (Of course remembering such an extra factor for the
product expression.)
\smallskip

Using algorithm~\ref{heightbase} for $\SL_n(\Z)$, we then can write the 
$\{1,\ldots,n\}\times\{1,\ldots,n\}$ minor $M$ of $e$ (respectively $e'$) as a
product of elementary matrices in dimension $n$.

As the generating set $\mvec{\tilde g}$ also contains
(product expressions for) matrices that act on this minor $M$ as elementary
matrices:
Let $w$ be an $\SL_{n}(\Z)$ word for $M$.
Evaluating $w$ in the generators $t_{i,j}t_{n+j,n+i}^{-1}\in\mvec{\tilde g}$
then gives an expression for $e$ in generators for $\Sp_{2n}(\Z)$.
\smallskip

It thus is sufficient to map an element
$e\in\Sp_{2n}(\Z)$ into $R$.
\medskip

An element $a=(a_{i,j})\in\Sp_{2n}(\Z)$ lies in $R$, if the height function
\[ h(a)=\sum_{i=1}^n \sum_{j=1}^n \left(a_{i,n+j}^2+a_{n+i,j}^2\right) \]
has value zero. This suggests that we can transform $e\in\Sp_{2n}(\Z)$ into
an element of $R$ by running algorithm~\ref{heightbase} with this new height
function $h$ (even though it does not have a unique minimal element).
\smallskip 

We thus modify the height-based approach of algorithm~\ref{heightbase} as
follows:
\begin{enumerate}
\item The stopping condition, in the outermost while-loop, is for $h(A)=0$
rather than $a=I$;
\item It returns not only the word expression, but also the reduced element
$a$;
\item the case $m\ge\|a\|$ first uses the above modification of the
algorithm that first tries an extended generating set, formed by adding
short products in the generators, before triggering an error if this also
found no reduction.
\end{enumerate}

Again, experiments with this approach failed, producing matrices that had
only a few nonzero entries in the top right and bottom left quadrant with no
way to also zero out these remaining entries. The goal to reduce all entries
at the same time led to a local, not global, minimum from which escape was
not possible.

To avoid such a behavior, we switched to a more localized reduction.
Based on the observation that single nonzero entries are hard to
clean out if the rest of their row is zero, we switch to an iterated
process, reducing row-by-row.
That is we define a series of
height functions by
\begin{eqnarray*}
h_0&=&0\\
h_i&=&h_{i-1}+\sum_{j=n+1}^{2n}a_{i,j}^2,\qquad\mbox{if $i\le n$}\\
h_i&=&h_{i-1}+\sum_{j=1}^{n}a_{i,j}^2,\qquad\mbox{if $i>n$}.
\end{eqnarray*}
We then run algorithm~\ref{heightbase} to reduce by height function $h_1$.
Afterwards, we reduce further with height function $h_2$ and so on, up to height function
$h_{2n}=h$. If the resulting matrix lies in $S$, we proceed as described
above, producing a word expression for $e$.
We call this approach Algorithm~5.
\medskip

We note that this improved heuristic succeeded in all examples we
tried (i.e found a matrix $a\in S$). We did not encounter a single example
in which this approach failed.

It also produced words of acceptable length.
Alas, proving these statements as general facts seems to be beyond the
capabilities of the author.
\medskip

What seems to be happening is that the localized heights are willing to accept
reduction step that reduce the current row, even if they grow the
entries in other places that are not covered by the height.

Contrary to the overall height function $h$, this approach thus does not
forbid a reduction to zero (which might produce a very small height
reduction), just because it combines with a growth of larger, not yet
reduced, entries of the matrix (note that an entry change $m$ to $m+1$
increases the height by roughly $(m+1)^2-m^2=2m+1$, while a reduction $1$ to
$0$ reduces by $1$ only).
\medskip

\section{Examples}
\label{xample}
As mentioned in the introduction, our main interest has been to obtain short
words for $\Sp$. We thus did not measure run times (which can be heavily
biased by setup costs or cleverness in avoiding duplicate calculations of
elements) systematically, but rather the quality of words obtained. This was
done in a {\sf GAP}~\cite{GAP4} implementation of the algorithms described
here, that is part of the author's routines for arithmetic groups, available
at~\url{www.math.colostate.edu/~hulpke/arithmetic.g}.
\smallskip

For a small example, section~7.2 in~\cite{longreid11}, using {\sf
Mathematica}, computes word expressions for selected elements of
$\SL_3(\Z)$, namely $X_0$ of length 8 and $Y_0$ of length 14;
algorithm~\ref{heightbase} obtained word expressions of length 7 and 13,
respectively.  Similarly an element $X_{-2}$ is given by a word of length 13
and $Z_{-2}$ by a word of length 16; algorithm~\ref{heightbase} calculated
expressions of lengths 16 and 10 respectively. The new approach thus
performs on par with an existing method.
\medskip

The next example is the group $G(3,4)$ from~\cite{hofmanvanstraten},
already considered in~\cite{detinkoflanneryhulpke17}. Using the
implementation in {\sf GAP} we construct a homomorphism from a finitely
presented version of $\Sp_4(\Z)$ to a matrix version, using the extended
generating set based on~\cite{birman71}. We also form $G(3,4)$ as a
matrix group. We then express (this uses the symplectic method) the
group generators as words, and form the subgroup $S$ of the finitely
presented $\Sp_4(\Z)$ that is generated by these words. We finally determine
the index $[\Sp_4(\Z):S]$ through a coset enumeration. (This calculation,
incidentally, independently verifies that $G(3,4)$ is arithmetic.)
\begin{verbatim}
gap> hom:=SPNZFP(4);
[ Y1, Y2, U1, U2, Z1 ] ->
[[[1,0,-1,0],[0,1,0,0],[0,0,1,0],[0,0,0,1]], [...]
gap> G34:=HofmannStraatenExample(3,4);
<matrix group with 2 generators>
gap> w:=List(GeneratorsOfGroup(G34),
> x->PreImagesRepresentative(hom,x));
[ U1*(U2^-1*U1*U2^-1)^2*Y1^-1*Y2^-1*U2^-1*Y2^-1
  *Z1*U2*Y2, Y2^-1 ]
gap> S:=Subgroup(Source(hom),w);;
gap> Index(Source(hom),S);
3110400
\end{verbatim}
In this example,  finding the word expressions (of length 14, respectively
1) takes $0.1$ seconds (while the coset enumeration confirming the index
takes about 4 minutes).

With algorithm~\ref{brute}, we verified
(in 20 minutes) that there is an expression for the first generator
of length $12$. Using this shorter word did not seem to have a
meaningful impact on the time required by the coset enumeration.
\bigskip

The input to all other experiments were matrices obtained as random words,
of a preselected length $len$, in the matrix generators of $\SL$,
respectively $\Sp$.  This produced matrices in the respective group for
which an upper bound for the length of a word expression was known. The
dimensions considered were chosen for be $\le 8$, as the motivating examples
from~\cite{detinkoflanneryhulpke17} do not exceed this bound.

We then calculated for each of the matrices a word expression, using the
algorithms described in this paper. If an algorithm produced a word of
length $a$ for a chosen input length $len$, we use the scaled ratio
$q:=100\cdot a/len$ as a a measure for the quality of the word expression
obtained. The diagrams given indicate a distribution of how often (the
ordinate) certain ratios $q$ (the abscissa) occur. (Incidentally, the
required runtime is reasonably approximated by this ratio,
as the fundamental step in all algorithms is to divide off one generator
matrix, building the word in steps of length one.)

The lengths considered were $20$ and larger which led to matrices whose
entries were frequently in the thousands or more. We therefore did not
attempt comparisons with algorithms~\ref{brute} or~\ref{congrimg}.

\begin{figure}
\includegraphics[width=4cm]{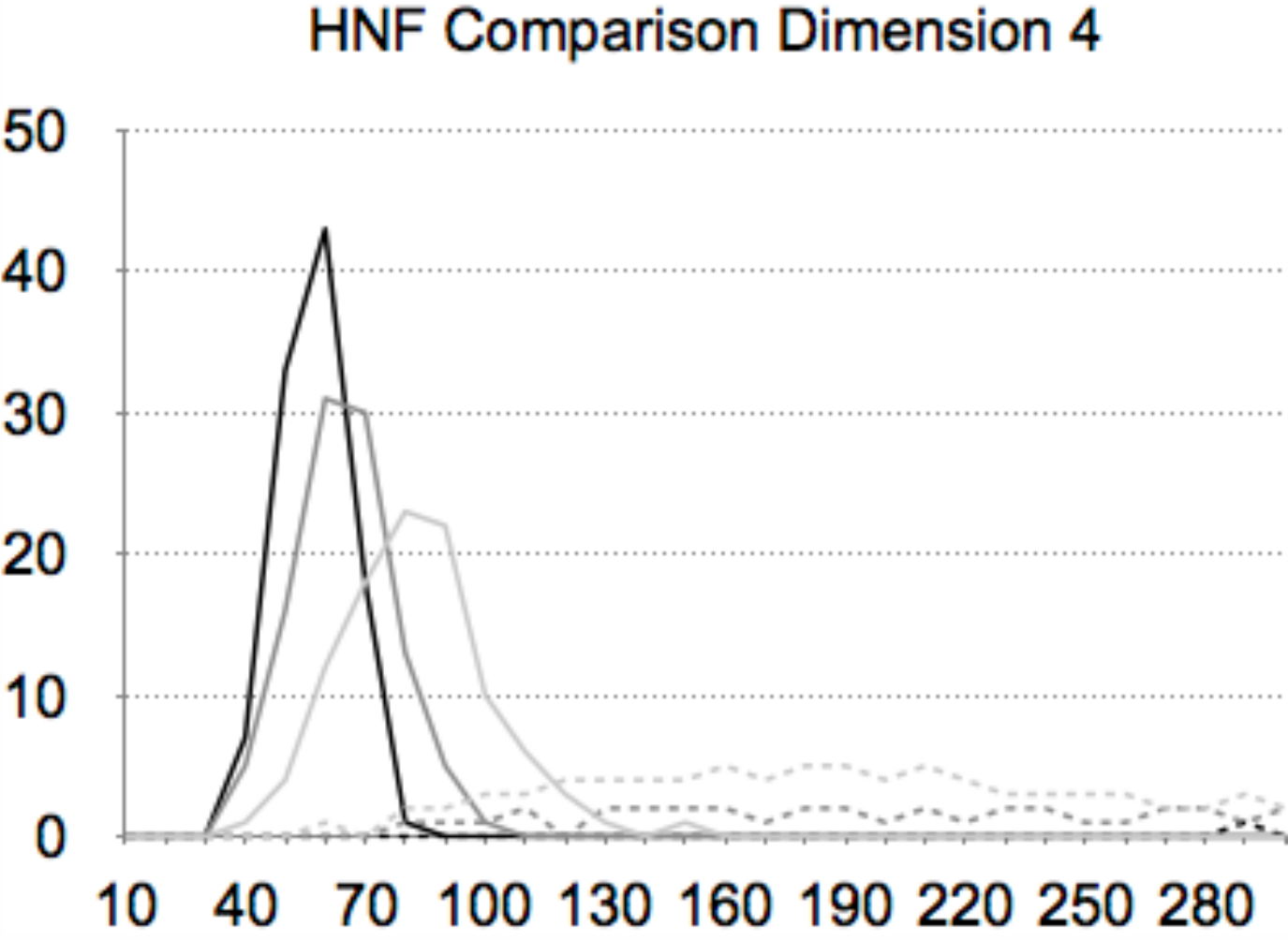}\quad
\includegraphics[width=4cm]{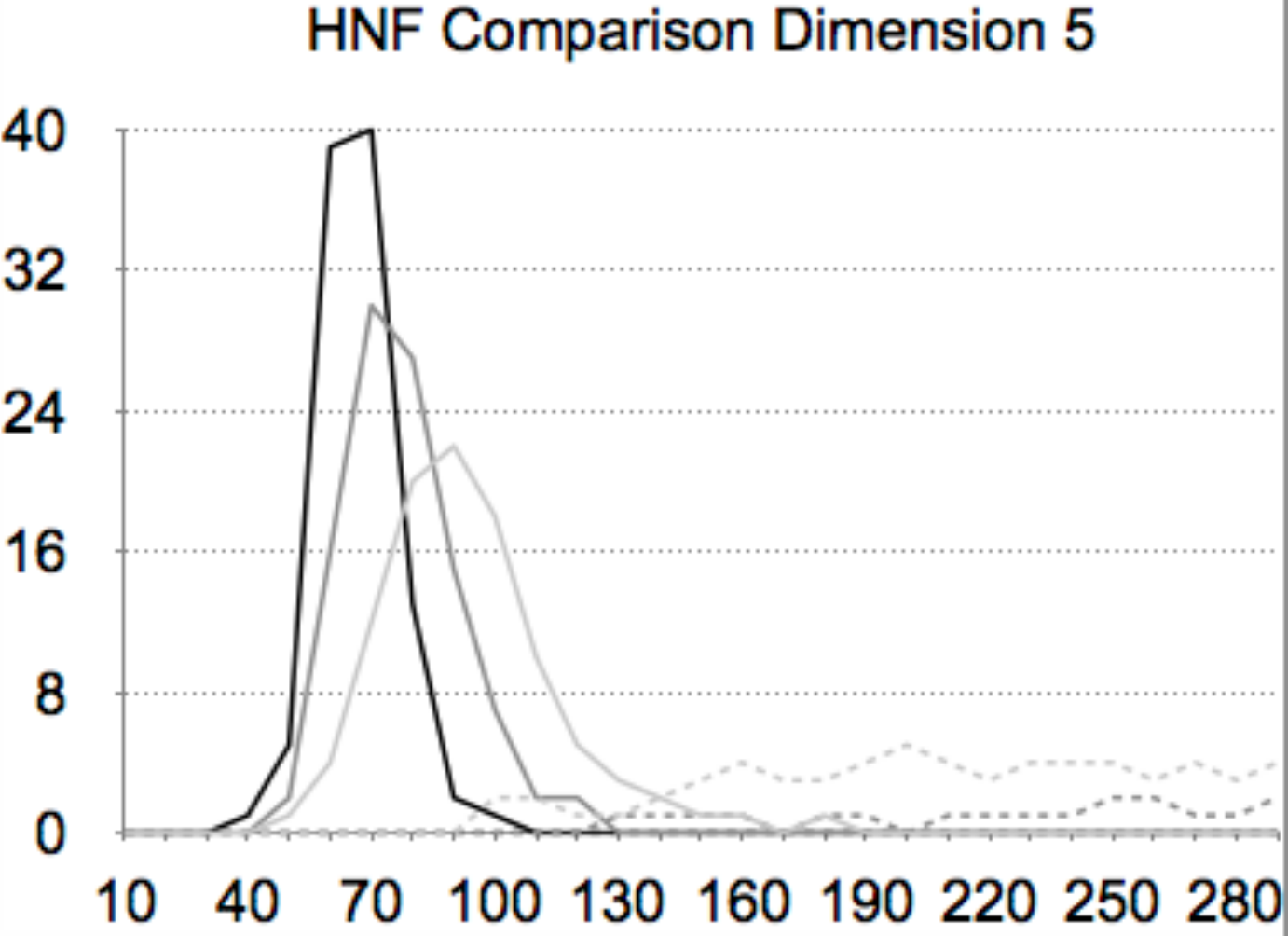}\\
\includegraphics[width=4cm]{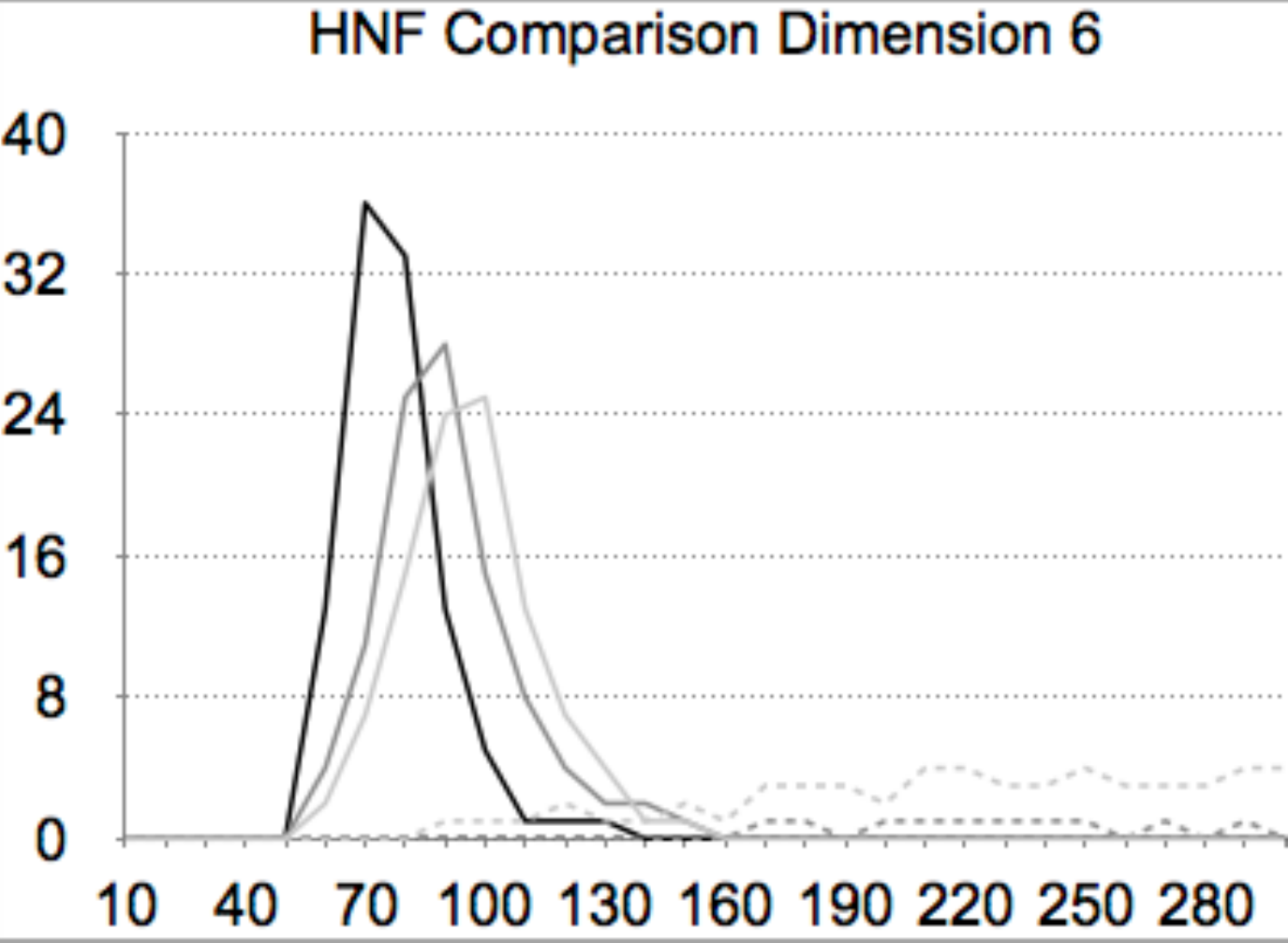}\quad
\includegraphics[width=4cm]{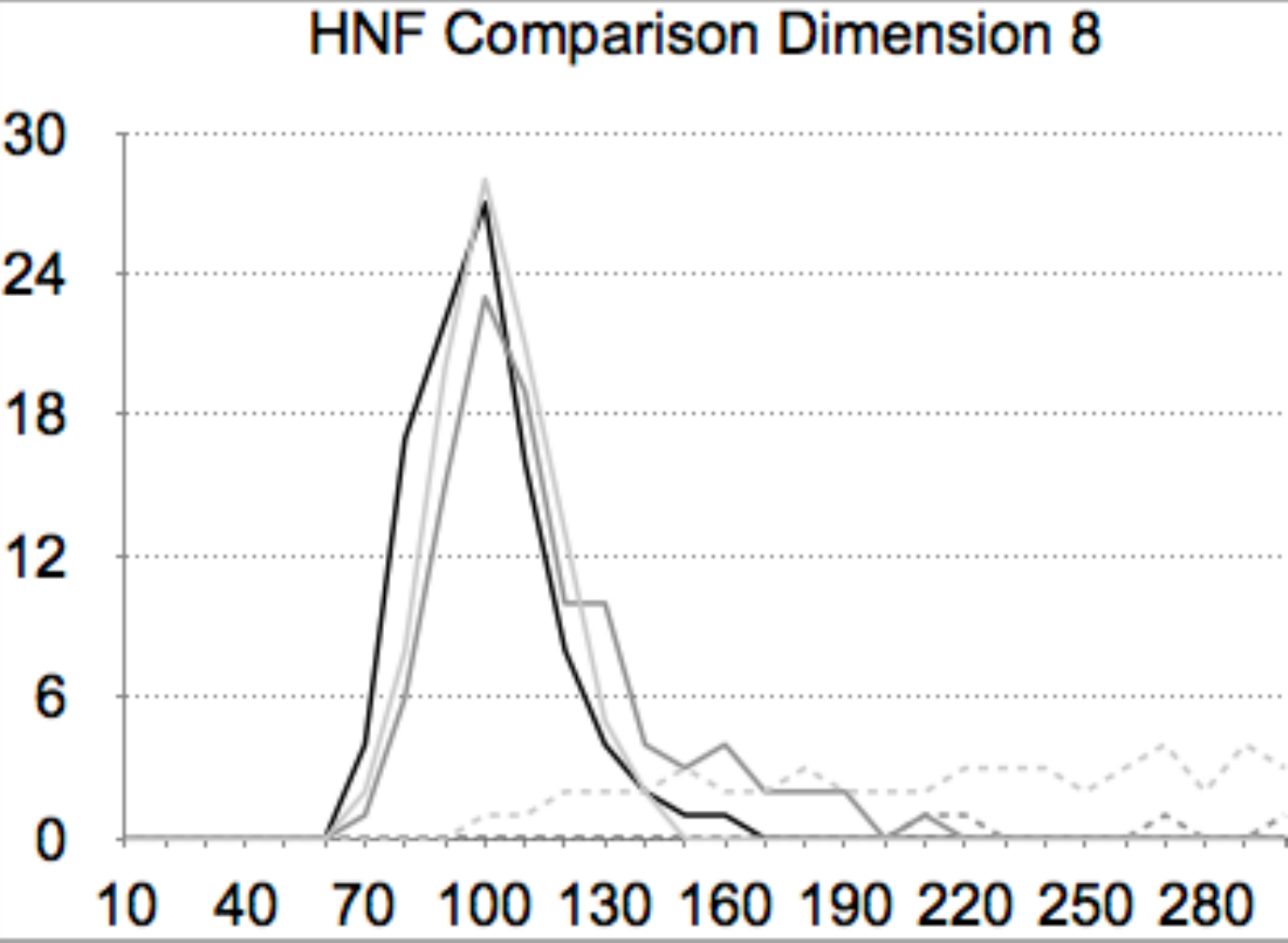}\\
\caption{Comparison between HNF-based algorithm~\ref{hnfbased} and
height-based algorithm~\ref{heightbase} for $\SL$}
\label{fighnf}
\end{figure}

As we only had time for a limited number of trials --- we used
$20000/len$ matrices of input length $len$ --- we discretized
the distribution in the following way to produce diagrams that  are easily
reproduced in print.
We grouped the
ratios $q$ into intervals of length $10$ each, and for each interval
calculated
the percentage of cases within the experiments for which the obtained ratio
fell into this interval. To allow for multiple experiments within one
diagram we did plot these results as piecewise linear curves (that somewhat
approximate a Gaussian distribution), rather than as bar graphs.
So for example in the top-left diagram in figure~\ref{fighnf} the continuous
black line indicates that about $5\%$ of experiments resulted in a ratio in
the interval $[40,50)$, $32\%$ in the interval $[50,60)$, $42\%$ in the
interval $[60,70)$ and (this is hard to see) $17\%$ in the interval
$[70,80)$, with the remaining $4\%$ of experiments resulting in ratios
outside this range (and too low to really show up in the diagram).
\medskip

The first series of experiments, given in figure~\ref{fighnf}, compares the
HNF-based algorithm~\ref{hnfbased} (dashed lines) with the height-based
algorithm~\ref{heightbase}, including the improvements by short products, on
matrices in $\SL_n(\Z)$ (continuous lines). We tested input lengths $20$,
$50$ and $100$ with darker colors representing longer input lengths, that is
$len=100$ is black, $len=50$ is mid-gray and $len=20$ is light gray. (For a
given input length. The {\em same set} of matrices was used for both algorithms.)

One immediately notes from the figures that the height based algorithm
produced results that (with some goodwill) can be
considered as approximations of a Gaussian distributions, centered not too
far off $100$.

The pure HNF-based algorithm instead produced a much wider spectrum of
results (the curves continue beyond the right edge of the diagram, which is
the reason the dashed black curves are practically invisible), with the average
length ratio becoming worse with longer word lengths.  Concretely, in the
case of dimension 4 and input length 100, the HNF-based algorithm produced
words whose length ranged between $290$ and $8,600,000$ (with an average of
$350,000$), making them useless in practice.

We thus
conclude (somewhat unsurprisingly, given what is known about integer normal form
calculations) that, at similar runtime, algorithm~\ref{heightbase} produces
significantly shorter words than a systematic HNF calculation.
\medskip

\begin{figure}
\includegraphics[width=4cm]{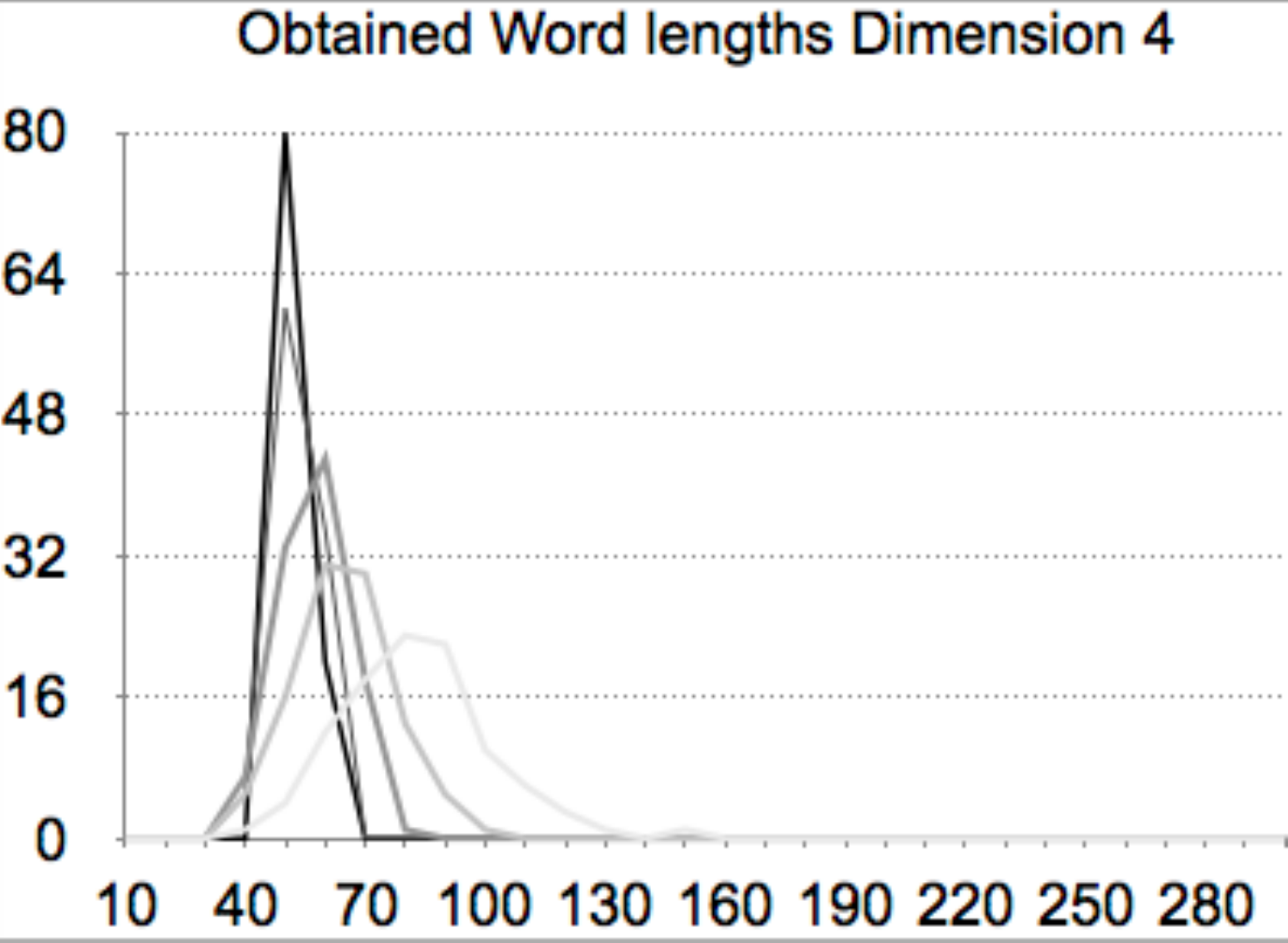}\quad
\includegraphics[width=4cm]{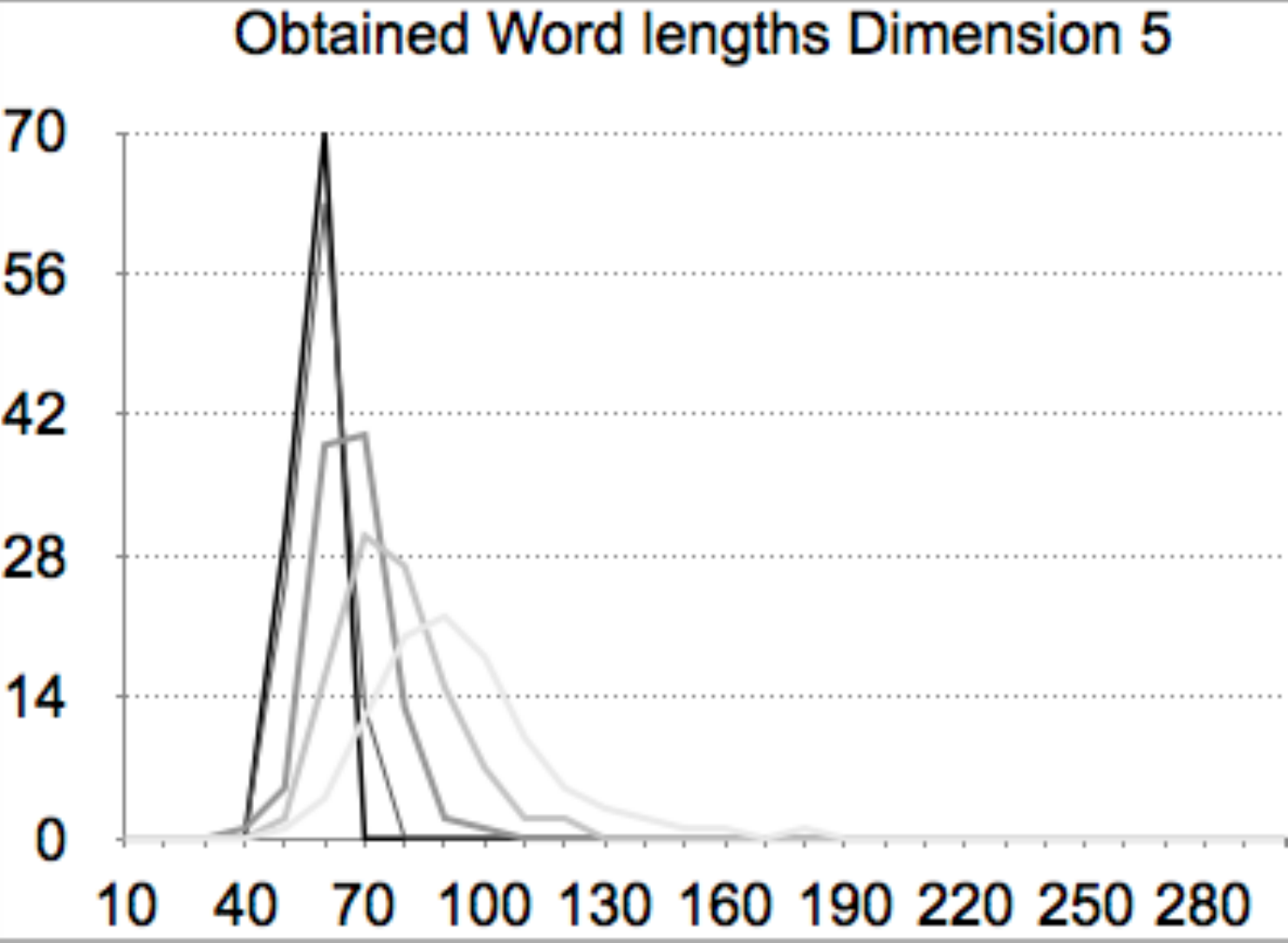}\\
\includegraphics[width=4cm]{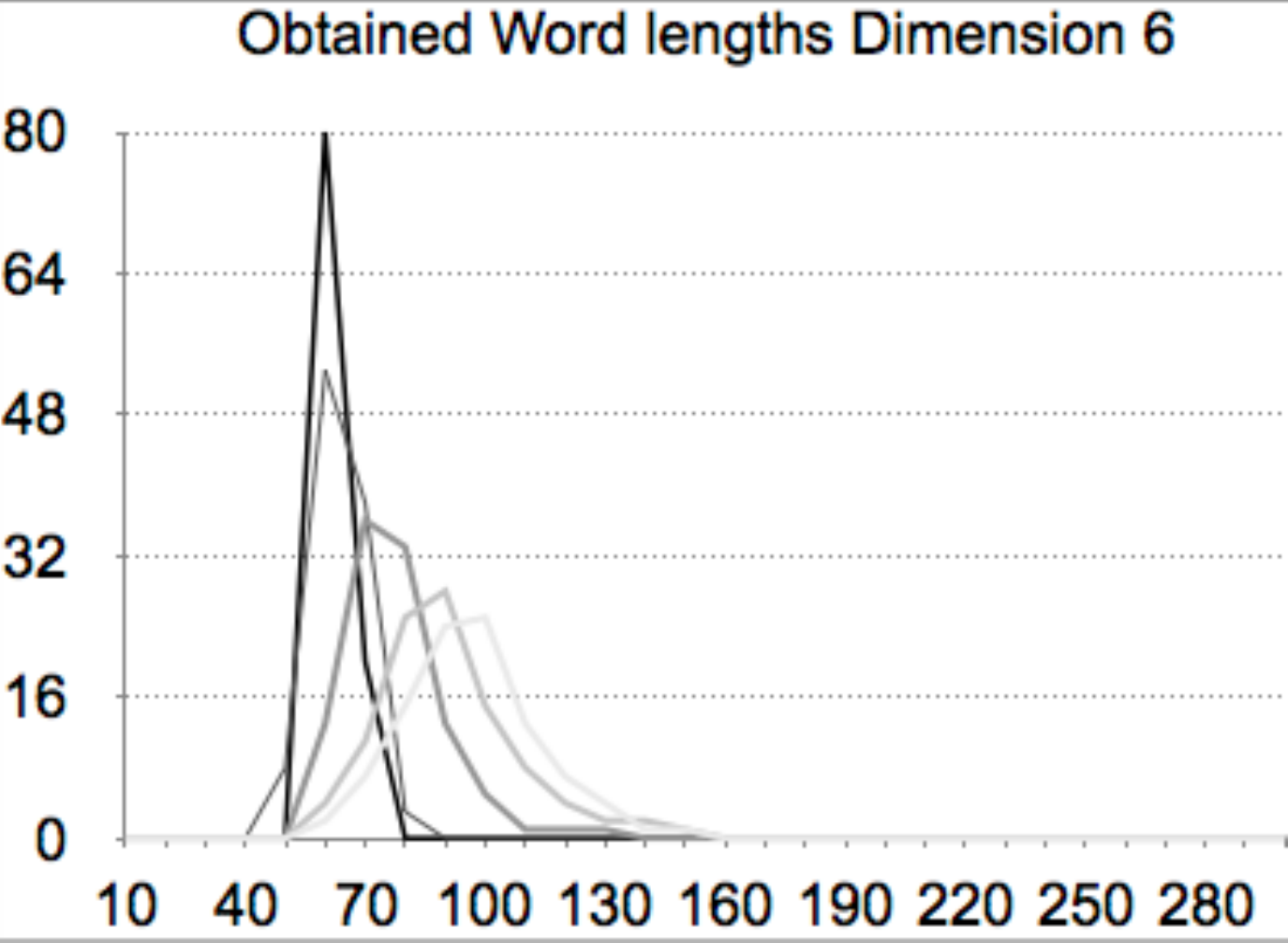}\quad
\includegraphics[width=4cm]{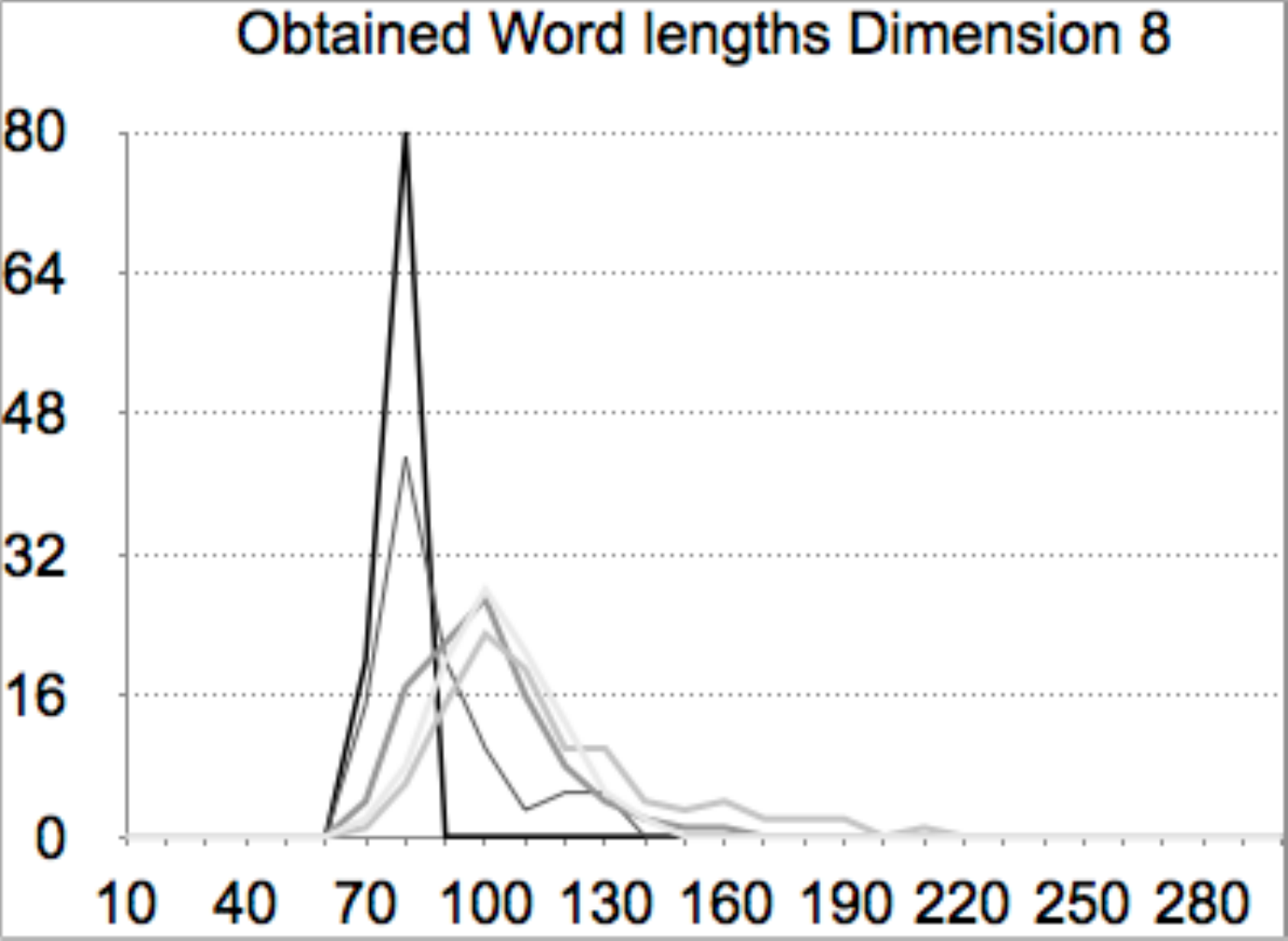}\\
\caption{Comparison of obtained word length for $\SL_n$ in different
dimensions}
\label{figsldim}
\end{figure}

In the second series of experiments in figure~\ref{figsldim} we considered only
algorithm~\ref{heightbase}, but for a broader set of lengths.
The input consisted of matrices given by random words of lengths $20$, $50$,
$100$, $500$ and $2000$, with darker colors again representing longer lengths.
Again we used $20000/len$ words of length $len$.

As before we observe curves that approximate a Gaussian distribution, with
peaks shifting to the left as the input length increases, and shifting to the
right as dimension increases. In the dimension range tested, both changes
are small enough to be considered as linear with small constant.

We did not try larger dimensions systematically, as calculations quickly
became unreasonably costly.
\smallskip

We note that the seeming improvement in the  resulting word lengths for
longer input might instead indicate that random words are less likely to be
optimally reduced as the length increases.
\bigskip

As for comparison with the optimal word length, this optimal length alas is
unknown in the examples (and because of memory limitations cannot be
determined for examples with input length $100$). Considering the rapid
growth observed for small lengths in the number of different elements that
can be expressed as words of a particular length, however it seems plausible
to have optimal length of the elements considered would differ from the
input length by a factor that is logarithmic in the word length rather than
linear.
\bigskip

The third series of experiments concerns elements of $\Sp$ for various
dimensions and lengths, using algorithm~5.

\begin{figure}
\includegraphics[width=4cm]{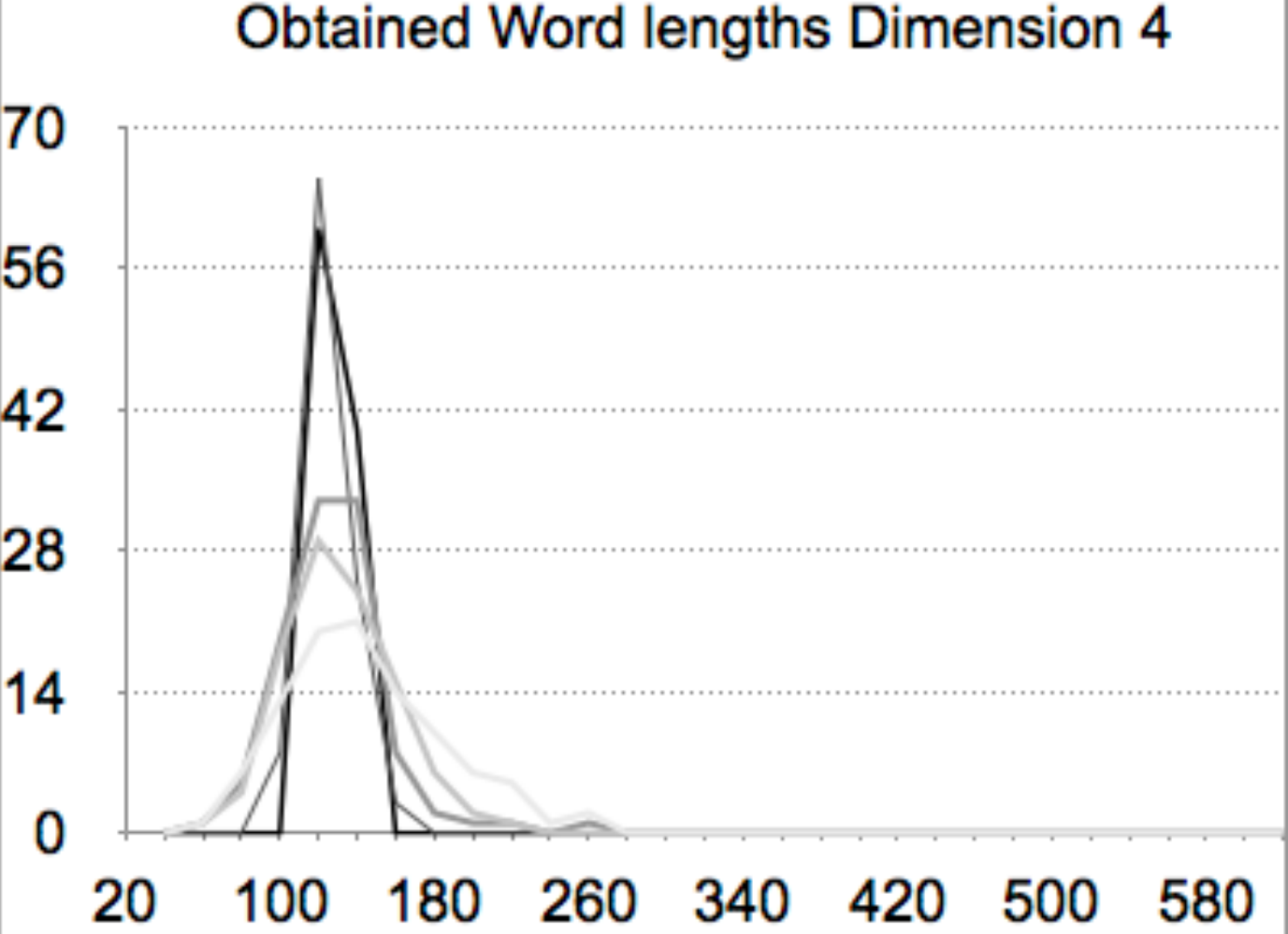}\quad
\includegraphics[width=4cm]{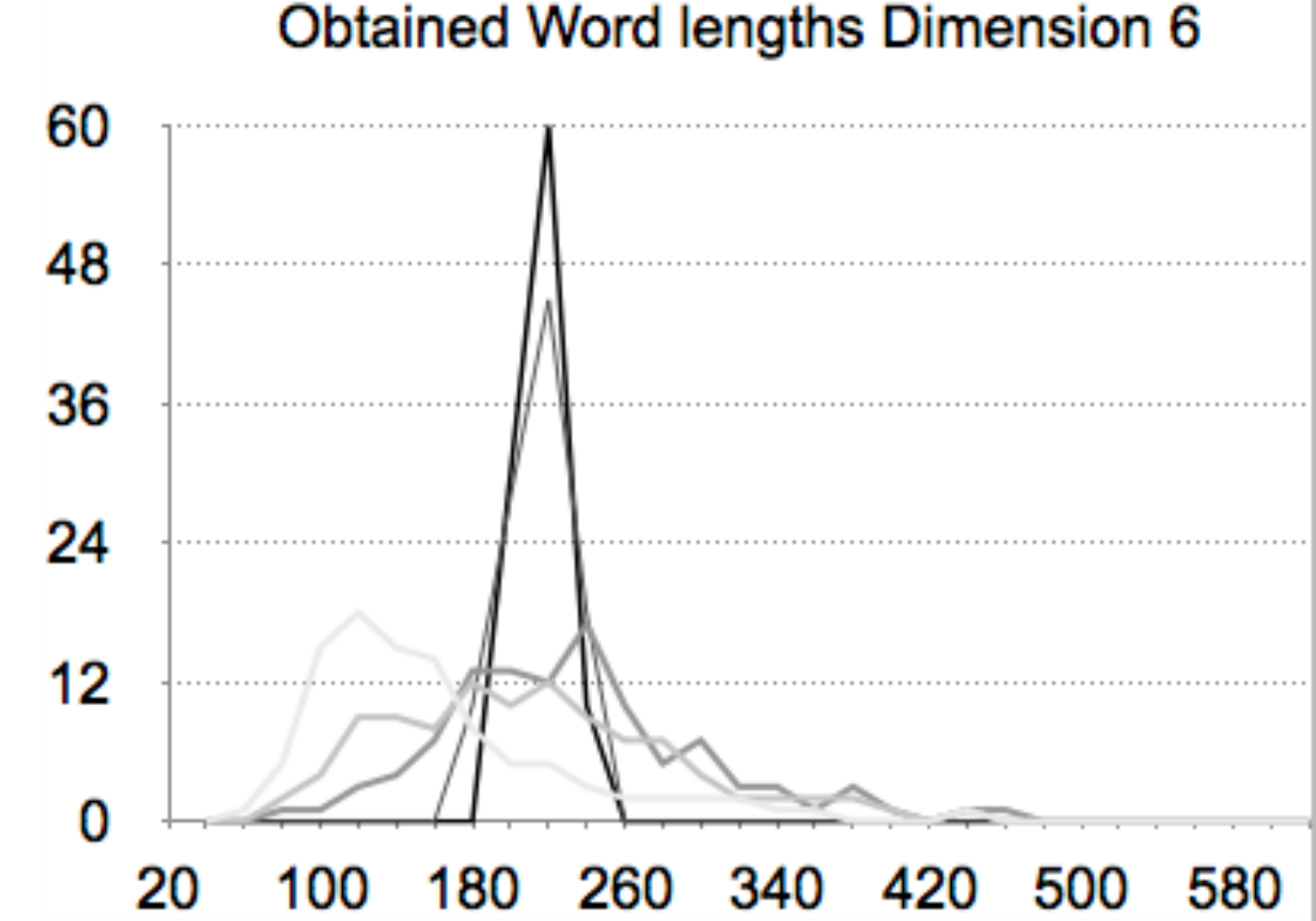}\\
\includegraphics[width=4cm]{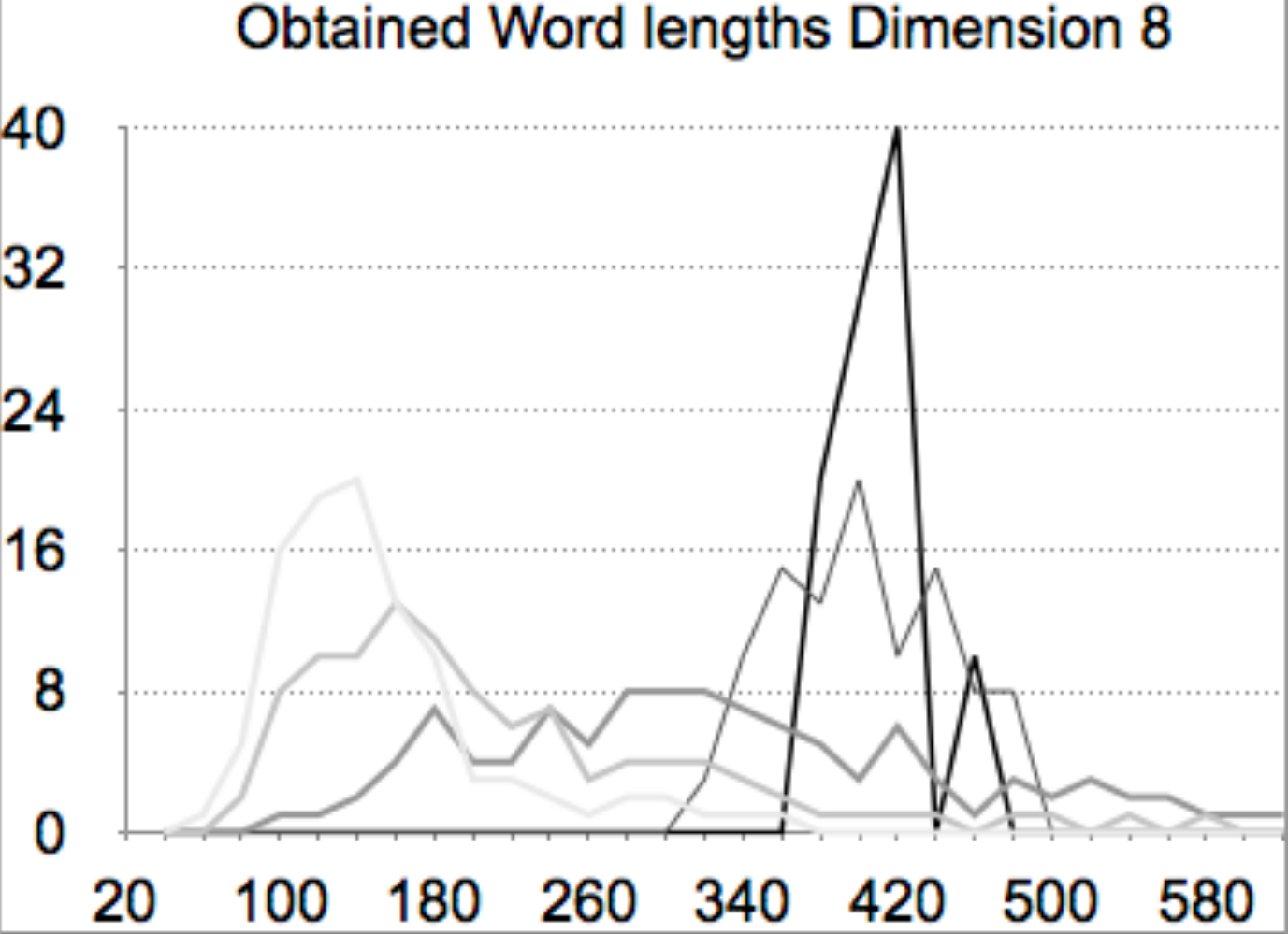}\\
\caption{Comparison of obtained word length for $\Sp_n$ in different
dimensions}
\label{figspdim}
\end{figure}

Figure~\ref{figspdim} gives the results of these experiments.
(Input lengths used and colors are as in the second
series.) For each of the random example matrices tested,
the approach found a factorization.

In dimension 4 the result is very similar as for $\SL$. With growing
dimension the behavior changes:
The larger number of generators acting locally on matrices make it more
likely that randomly chosen generators commute. If the word length is short
one can almost read off the generators involved from the positions of
nonzero matrix entries.

Longer words in larger dimensions however show an increased widening of the
bell shape and a shift of the peak towards significantly longer words --
about $200\%$ for dimension $6$ and $400\%$ for dimension 8.  This seems to
indicate that the approach is feasible for small dimensions (not least for
the lack of alternatives) with word lengths not increasing by too much and
no observed failure, but that for larger dimensions the ratio to optimal
word length gets exponentially worse.

\section{Closing remarks}

We have seen practically feasible methods to express elements of $\SL_n(\Z)$
and $\Sp_{2n}(\Z)$ (in small dimensions) as products in particular
(standard)
generating sets. 

What is 
clearly lacking is a proof (and not just experimental evidence) of the
approach succeeding in general for $\Sp$, as well as of the 
produced words being not too worse than the minimal word lengths for the
matrices.
(The latter seems difficult as the algorithm for $\SL$ which is proven to
succeed within limited memory -- the HNF-based one -- produced words of
unusable length.) Even without such a proof the heuristic presented will be
useful, as long as it produces a result.

The tools motivating our approach were taken from integral matrix normal
forms. This raises the question on whether further synergies in either way
can be obtained from these problems. A first caveat is that the normal form
in the factorization case is always the identity matrix, and that any
experiments done here were in tiny dimensions compared with those usually
considered for normal forms.

What might be more promising (but we have not investigated) is a relation
between word length and size of matrix entries for the transforming
matrices for e.g. the Smith Normal Form.
We observed that an initial norm-based global reduction of matrix
norms produced significantly shorter words. If this can be translated to
smaller matrix entries, it would be useful for applications such as the
homomorphisms to abelianizations $G/G'$ of finitely presented groups.

\section{Acknowledgments}
The author's work has been supported in part
by Simons Foundation Collaboration Grants~244502 and~524518 which are gratefully
acknowledged. The author also would like to thank the anonymous referees for
their helpful remarks.

\bibliographystyle{plain}

\end{document}